\documentclass[reqno,centertags,draft,10pt]{amsart}
\usepackage{amssymb}
\usepackage{upref}
\newcommand{\bbN}{{\mathbb{N}}}
\newcommand{\bbR}{{\mathbb{R}}}

\newcommand{\bbC}{{\mathbb{C}}}

\newcommand{\calC}{{\mathcal C}}

\newcommand{\calD}{{\mathcal D}}

\newcommand{\dott}{\,\cdot\,}
\newcommand{\no}{\nonumber}
\newcommand{\lb}{\label}
\newcommand{\f}{\frac}

\newcommand{\ti}{\tilde  }

\newcommand{\loc}{\text{\rm{loc}}}
\newcommand{\dom}{\text{\rm{dom}}}
\newcommand{\rank}{\text{\rm{rank}}}
\newcommand{\ess}{\text{\rm{ess}}}
\newcommand{\ac}{\text{\rm{ac}}}
\newcommand{\singc}{\text{\rm{sc}}}
\newcommand{\pp}{\text{\rm{pp}}}
\newcommand{\supp}{\text{\rm{supp}}}
\newcommand{\AC}{\text{\rm{AC}}}
\newcommand{\bi}{\bibitem}

\renewcommand{\Re}{\text{\rm Re}}
\renewcommand{\Im}{\text{\rm Im}}

\allowdisplaybreaks
\numberwithin{equation}{section}

\newtheorem{theorem}{Theorem}[section]

\newtheorem{corollary}[theorem]{Corollary}
\newtheorem{hypothesis}[theorem]{Hypothesis}
\theoremstyle{definition}
\newtheorem{definition}[theorem]{Definition}

\theoremstyle{remark}
\newtheorem{remark}[theorem]{Remark}

\newcommand{\abs}[1]{\lvert#1\rvert}

\begin{document}
\title[Borg-type theorems]{Borg-type theorems for matrix-valued
Schr\"{o}dinger operators}
\author[Clark, Gesztesy, Holden, and Levitan]{Steve Clark, 
Fritz Gesztesy, Helge Holden, \linebreak and Boris
M.~Levitan}
\address{Department of Mathematics and Statistics, University 
of Missouri-Rolla, Rolla, MO 65409, USA}
\email{sclark@umr.edu}
\urladdr{http://www.umr.edu/\~{}clark}
\address{Department of Mathematics,
University of Missouri,
Columbia, MO 65211, USA}
\email{fritz@math.missouri.edu}
\urladdr{http://www.math.missouri.edu/people/fgesztesy.html}
\address{Department of Mathematical Sciences, 
Norwegian University of
Science and Technology, N--7491 Trondheim, Norway}
\email{holden@math.ntnu.no}
\urladdr{http://www.math.ntnu.no/\~{}holden/}
\address{School of Mathematics, University of 
Minnesota, Minneapolis,
MN 55455, USA}
\email{levit009@maroon.tc.umn.edu}
\thanks{Research supported in part by the Norwegian 
Research Council.}
\date{\today}

\begin{abstract}
A Borg-type uniqueness theorem for matrix-valued 
Schr\"odinger
operators is proved. More precisely, assuming 
a reflectionless
potential matrix and spectrum
a half-line $[0,\infty)$, we derive triviality 
of the potential matrix. 
Our approach is based on trace formulas and 
matrix-valued
Herglotz representation theorems. As a by-product 
of our techniques, we
obtain an extension of Borg's classical result 
from the class of periodic
scalar potentials to the class of reflectionless 
matrix-valued
potentials.
\end{abstract}

\maketitle
\section{Introduction}\lb{s1}
The principal aim of this paper is to advocate 
a new proof of
Borg's \cite{Bo46}
uniqueness result for periodic one-dimensional 
Schr\"odinger operators in
$L^2(\bbR)$ and to extend it to general 
matrix-valued Schr\"odinger 
operators in $L^2(\bbR)^m$. In order to 
describe Borg's result and our
generalizations of it (see 
Section~\ref{s4} for more
details), we need a few preparations. Let 
$H^{p,q}(\bbR)$, $p,q \in
\bbN$ denote the standard Sobolev spaces and 
$\sigma(\dott)$ abbreviate
the spectrum. Assuming
\begin{align}
&q \in \AC_{\loc} (\bbR) \text{ to be 
real-valued,} \lb{1.1} \\
& {\sigma}(h)=[e_0,\infty) \text{ for some 
$e_0\in \bbR$,} \lb{1.2} \\
&q \text{ periodic}, \lb{1.3}
\end{align}
with $h$ on $H^{2,2}(\bbR)$ denoting the usual
self-adjoint realization of the differential 
expression
$-\f{d^2}{dx^2}+q(x)$ in $L^2(\bbR)$, 
Borg \cite{Bo46} proved the
uniqueness result
\begin{equation}
q(x)=e_0 \text{ for all $x \in \bbR$.} \lb{1.4}
\end{equation}
(Actually, Borg only assumed 
$q \in L^2_{\loc} (\bbR)$ and hence
obtained $q=e_0 \text{ a.e.}$ in \eqref{1.4} but we 
will temporarily
ignore this for simplicity of exposition.)

Next, consider matrix-valued Schr\"odinger operators 
$H$ in
$L^2(\bbR)^m$, $m\in \bbN$ associated with 
differential expressions
\begin{equation}
-I_m\f{d^2}{dx^2} + Q(x), \quad x \in \bbR, \lb{1.5}
\end{equation}
where $I_m$ is the identity matrix in $\bbC^m$, $Q \in
\AC_{\loc}(\bbR)^{m \times m}$, and $Q(x)=Q(x)^{*}$ 
for all $x\in
\bbR.$

In Section~\ref{s4} we will prove a uniqueness result 
of the type
\eqref{1.4} for the matrix-valued operators $H.$ 
Our main tool will
be a trace formula of the type
\begin{equation}
Q(x)=E_0 I_m+\lim_{z\to i\infty}\int_{E_0}^\infty 
d\lambda\,
z^2(\lambda-z)^{-2}(I_m-2\Xi(\lambda,x)), 
\quad x\in\bbR \lb{1.7}
\end{equation}
for Schr\"odinger operators. Here 
$\Xi(\lambda,x)$ is a
self-adjoint $m \times m$ matrix satisfying
\begin{equation}
0 \leq \Xi(\lambda,x) \leq I_m \lb{1.8}
\end{equation}
and $E_0$ denotes the infimum of the spectrum of $H$. 
Given the trace
formula \eqref{1.7}, our idea of extending Borg's 
uniqueness theorem
to matrix-valued Schr\"odinger operators now 
becomes very 
simple. Suppose
the analog of condition \eqref{1.2}, that is,
\begin{equation}
\sigma(H)=[E_0,\infty) \text{ for some } E_0 
\in \bbR \lb{1.9}
\end{equation}
and instead of the periodicity condition \eqref{1.3} 
assume that for
all $x\in\bbR$,
\begin{equation}
\Xi (\lambda,x)= \f{1}{2} I_m \text{ for a.e. } 
\lambda \in
[E_0,\infty).
\lb{1.10}
\end{equation}
Then the trace formula \eqref{1.7} immediately 
yields
\begin{equation}
Q(x)=E_0 I_m \lb{1.11}
\end{equation}
and hence a desired generalization of Borg's 
result \eqref{1.4}. We
will show in Section~\ref{s4} that periodicity of 
$Q(x)$, together
with the assumption that $H$ has uniform (maximal) 
spectral
multiplicity $2m$, indeed implies 
condition \eqref{1.10}. This
recovers a generalization of Borg's theorem to periodic
matrix-valued Schr\"odinger operators by 
Depr\'es \cite{De95}, which
partly motivated our present work. Consequently, 
assumption
\eqref{1.10} is a
proper extension of the periodicity 
requirement \eqref{1.3}.

More generally, if for all $x\in\bbR$,
\begin{equation}
\Xi(\lambda,x)=\f12 I_m \text{ for a.e. } \lambda \in
\sigma_{\ess}(H) \lb{1.12}
\end{equation}
($\sigma_{\ess} (\dott )$ denoting the essential 
spectrum), we 
shall call
$Q(x)$ a {\it reflectionless} potential following the 
traditional
terminology in the scalar case $m=1$ (cf. the 
discussion in
Section~\ref{s4}). Thus
reflectionless potentials $Q(x)$, in connection with 
the spectral
assumption \eqref{1.9}, are the prime candidates for 
Borg-type
theorems. 

Finally, we briefly sketch the content of each section. 
Section~\ref{s2} provides the basic background results on 
matrix-valued
Schr\"odinger operators. Following a series of
papers by Hinton and Shaw \cite{HS81}--\cite{HS84}, \cite{HS86}, 
and with a view toward a future treatment of Dirac-type 
operators, we treat Schr\"odinger operators as special 
Hamiltonian systems and briefly recall the corresponding
Weyl-Titchmarsh and spectral theory.
In Section~\ref{s3} we derive new trace formulas for
matrix-valued Schr\"odinger operators using
appropriate Herglotz representation results for 
a diagonal Green's
matrix as discussed in Section~\ref{s2}. In our principal 
Section~\ref{s4} we finally derive the extension of 
Borg-type 
theorems to
matrix-valued Schr\"odinger operators. We also
provide a criterion for a potential to be 
reflectionless and
 close with an application to the case of periodic 
potentials.

Our results can be viewed as a first (and rather 
modest) step
toward the construction of isospectral manifolds of 
certain
classes of matrix-valued potentials for
Schr\"odinger operators. Especially, one might
think of the class of periodic (possibly reflectionless)
potentials, see \cite{Ca98a}, \cite{Ca99}. 
Moreover, our results are relevant in
the context of matrix-valued hierarchies of integrable 
evolution equations
(i.e., soliton equations) and we refer the reader 
to \cite{AK90},
\cite{Ch96}, \cite{Di91}, \cite{Di97}, \cite{DK97}, 
\cite{Du83}, \cite{GD77},
\cite{Ma78}, \cite{MO82}, \cite{OMG81}, \cite{WK74},
and the vast literature therein. For related work on 
trace formulas, spectral properties of matrix-valued 
Schr\"odinger 
operators, and uniqueness theorems see, for 
instance, \cite{AK92}, \cite{Ca98}, 
\cite{Ca98a}, \cite{Ca99}, \cite{Ch99}, \cite{Ch99a}, 
\cite{CS97}, \cite{JL98}, \cite{JL99}, \cite{Kh77}, 
\cite{MV81}, \cite{MV89}, \cite{Ma94}, \cite{Ma98}, 
\cite{Pa95}, \cite{Ro63}, \cite{Ry99}, \cite{Ry99a}, 
\cite{Sa94}, \cite{Ve83}, \cite{We87}.

The present paper focuses on matrix-valued
Schr\"odinger operators; corrersponding extensions 
to Dirac-type operators 
will appear elsewhere.

\section{Matrix-Valued Schr\"{o}dinger Operators} 
\lb{s2}

In this section we briefly recall the 
Weyl--Titchmarsh theory
for matrix-valued Schr\"{o}dinger operators. In 
view of a future treatment of Dirac operators we use 
a unifying approach representing Schr\"odinger 
operators as special cases of Hamiltonian
systems and
hence develop the theory from that point of view. 
Throughout 
this paper, all matrices will be considered over the 
field of 
complex numbers $\bbC.$

The basic assumption of this paper will be the following.

\begin{hypothesis}\lb{hyp2.1}
Fix $m\in\bbN$, $n=2m$, and define the $n\times n$ matrix
\begin{equation}
J=\begin{pmatrix}0& -I_m \\ I_m & 0  \end{pmatrix}. \lb{2.1}
\end{equation}
Suppose $Q=Q^*\in L_{\loc}^1(\bbR)^{m\times m}$ and
introduce the $n\times n$ matrices
\begin{equation}
A=\begin{pmatrix} I_m & 0 \\ 0& 0\end{pmatrix}, \quad
B(x)=\begin{pmatrix}-Q(x) & 0 \\0 & I_m\end{pmatrix}. 
\lb{2.2}
\end{equation}
\end{hypothesis}

Given Hypothesis~\ref{hyp2.1} we consider the 
Hamiltonian system
\begin{equation}
J\psi^\prime(z,x)=(zA+B(x))\psi(z,x), 
\quad z\in\bbC \lb{2.4}
\end{equation}
for a.e.\ $x\in\bbR$, where $z\in\bbC$ plays the 
role of a spectral
parameter and $\psi(z,x)$ is
assumed to satisfy
\begin{equation}
\psi(z,\dott) \in \AC_{\loc}(\bbR)^n. \lb{2.5}
\end{equation}
Here $I_p$ denotes the identity matrix in $\bbC^p$ for 
$p\in\bbN$, $M^*$
the adjoint (i.e., complex
conjugate transpose), $M^t$ the transpose of the 
matrix $M$, and
$\AC_{\loc}(\bbR)$ denotes the
set of locally absolutely continuous functions 
on $\bbR$.  
At times it will
be convenient to
consider  an $n\times r$ solution matrix 
in \eqref{2.4}, with
$r=1,\dots,n$, which will then be
denoted by $\Psi(z,x)$ and assumed to satisfy
$\Psi(z,\dott)\in\AC_{\loc}(\bbR)^{n\times r}$.

Hypothesis~\ref{hyp2.1} governs the case of
matrix-valued Schr\"{o}dinger
operators. In fact, equation \eqref{2.4} becomes 
equivalent to
\begin{align}
-\psi_1^{\prime\prime}(z,x)+Q(x)\psi_1(z,x)& =z\psi_1(z,x), 
\lb{2.6} \\
\psi_2(z,x)&=\psi_1^\prime(z,x), \lb{2.7}
\end{align}
where
\begin{equation}
\psi(z,x)=(\psi_1(z,x),\psi_2(z,x))^t. \lb{2.8}
\end{equation}
Here it is assumed that
\begin{equation}
\psi_j(z,\dott)\in\AC_{\loc}(\bbR)^m, \quad j=1,2. \lb{2.9}
\end{equation}

In order to recall the limit point and limit  circle cases 
associated with
the Hamiltonian
system \eqref{2.4}, we introduce the notation 
($ -\infty\le a< b \le \infty$)
\begin{align}
L_A^2((a,b))&=\{\phi:(a,b)\to\bbC^n  \mid \int_a^b dx\,
(\phi(x),A\phi(x))_{\bbC^n}<\infty \}, \lb{2.10} \\
N(z,\infty)&=
\{\phi\in L_A^2((c,\infty)) \mid J\phi^\prime=(zA+B)\phi
\text{ a.e. on $(c,\infty)$} \}, \lb{2.11}\\
N(z,-\infty)&=
\{\phi\in L_A^2((-\infty,c)) \mid J\phi^\prime=(zA+B)\phi
\text{ a.e. on $(-\infty,c)$} \}, \lb{2.11a}
\end{align}
for some $c\in\bbR$ and $z\in\bbC$. (Here
$(\phi,\psi)_{\bbC^n}=\sum_{j=1}^n \bar\phi_j\psi_j$
denotes the standard scalar product in $\bbC^n$, 
writing $\chi\in\bbC^n$ as 
$\chi=(\chi_1,\dots,\chi_n)^t,$ etc.)  Both 
dimensions of the spaces in
\eqref{2.11} and \eqref{2.11a},
$\dim_\bbC(N(z,\infty))$ and 
$\dim_\bbC(N(z,-\infty))$, are 
constant for
$z\in\bbC_\pm=\{\zeta\in\bbC
\mid \Im(\zeta)\gtrless 0 \}$, see for 
instance \cite{At64},
\cite{KR74}, who prove this fact for much more general 
Hamiltonian systems.
Hence one
defines  the Hamiltonian system \eqref{2.4} to 
be in the 
limit point (l.p.)
case at $\pm\infty$ if
\begin{equation}
\dim_\bbC(N(z,\pm\infty))=n/2 \text{  for all
$z\in\bbC\backslash\bbR$},  \lb{2.12}
\end{equation}
and in the limit circle case (l.c.) case at 
$\pm\infty$  if
\begin{equation}
\dim_\bbC(N(z,\pm\infty))=n \text{  for all $z\in\bbC$}.  
\lb{2.13}
\end{equation}
Later on we will introduce Schr\"{o}dinger operators 
$H$ in
$L^2(\bbR)^m$ associated with \eqref{2.4} and 
Hypothesis~\ref{hyp2.1} and see that 
the l.p.\ and l.c.\ notions
for the Hamiltonian system
\eqref{2.4} and the operators $H$ coincide. From this 
point on we
only consider the l.p.\
case at $\pm\infty$ and hence work with the 
following assumption.

\begin{hypothesis} \lb{hyp2.2}
In addition to Hypothesis~\ref{hyp2.1} suppose the 
Hamiltonian system
\eqref{2.4} to be in the l.p.\
case at $\pm\infty$.
\end{hypothesis}

Next we briefly turn to Weyl--Titchmarsh theory associated 
with \eqref{2.4}
and briefly recall some
of the results developed by Hinton and Shaw in a series of 
papers devoted
to spectral theory of
(singular) Hamiltonian systems \cite{HS81}--\cite{HS84}, 
\cite{HS86} (see
also \cite{Kr89a},
\cite{Kr89b}). While they discuss
\eqref{2.4} under much more general hypotheses on $A(x)$ and 
$B(x)$, we
here confine ourselves to
the special cases of  matrix-valued Schr\"{o}dinger systems
governed by Hypothesis~\ref{hyp2.2}.  Let
$\Psi(z,x,x_0)$ be a normalized fundamental system of solutions of
\eqref{2.4} at some
$x_0\in\bbR$, that is, $\Psi(z,x,x_0)$ satisfies
\begin{equation}
J\Psi^\prime(z,x)=(zA+B(x))\Psi(z,x), \quad z\in\bbC \lb{2.14}
\end{equation}
for a.e.\ $x\in\bbR$, and
\begin{equation}
\Psi(z,x_0,x_0)=I_n. \lb{2.15}
\end{equation}
Moreover, we partition $\Psi(z,x,x_0)$ as
\begin{equation}
\Psi(z,x,x_0)=(\psi_{j,k}(z,x,x_0))_{j,k=1}^2
=\begin{pmatrix}\theta_1(z,x,x_0) & \phi_1(z,x,x_0)\\
\theta_2(z,x,x_0)& \phi_2(z,x,x_0)\end{pmatrix}, \lb{2.16}
\end{equation}
where $\theta_j(z,x,x_0)$ and $\phi_j(z,x,x_0)$ for 
$j=1,2$ are 
$m\times m$
matrices, entire with
respect to $z\in\bbC$, and normalized according to 
\eqref{2.15}, 
that is,
\begin{equation}
\theta_1(z,x_0,x_0)=\phi_2(z,x_0,x_0)=I_m, \quad
\theta_2(z,x_0,x_0)=\phi_1(z,x_0,x_0)=0. \lb{2.17}
\end{equation}
(We recall that $\theta_2(z,x,x_0)=\theta_1^\prime (z,x,x_0)$ 
and $\phi_2(z,x,x_0)=\phi_1^\prime(z,x,x_0)$ in the present 
case of Schr\"odinger operators, cf.~\eqref{2.7}.)
One can prove \cite{HS81} that
\begin{equation}
\det(\phi_1(z,x,x_0))\neq 0 \text{ for 
$x\in\bbR\backslash\{x_0\}, 
\, z\in\bbC\backslash\bbR$} \lb{2.18}
\end{equation}
so that
\begin{equation}
M_{\pm,R}(z,x_0)=-\phi_1(z,R,x_0)^{-1}\theta_1(z,R,x_0), 
\quad R\gtrless
x_0, \,
z\in\bbC\backslash\bbR\lb{2.19}
\end{equation}
are well-defined.  Due to the assumption of the l.p.\ case 
at 
$\pm\infty$,
one obtains the
existence of the following limits \cite{HS81}, \cite{HS83}, 
\cite{HS84},
\cite{Or76}, \cite{Sa94a},
\begin{equation}
M_{\pm}(z,x_0)=\lim_{R\to\pm\infty} M_{\pm,R}(z,x_0), \quad
z\in\bbC\backslash\bbR.
\lb{2.20}
\end{equation}
$M_{+}(z,x_0)$ (resp.~$M_{-}(z,x_0)$) represent the 
half-line Weyl--Titchmarsh 
matrices
associated with \eqref{2.4} and
the interval $(x_0,\infty)$ (resp.~$(-\infty,x_0)$).  
For later reference we summarize the
principal results on
$M_{\pm}(z,x_0)$ in the following theorem.
\begin{theorem}[\cite{AD56}, \cite{Ca76}, \cite{GT97}, 
\cite{HS81},
\cite{HS82}, \cite{HS86},
\cite{KS88}]
\lb{thm2.3} Assume Hypothesis~\ref{hyp2.2} and let  
$z\in\bbC\backslash\bbR$, and
$x_0\in\bbR$.  Then \\
(i) $\pm M_{\pm}(z,x_0)$ is a matrix-valued Herglotz 
function of rank $m$.
In particular,
\begin{gather}
\Im(\pm M_{\pm}(z,x_0)) > 0, \quad z\in\bbC_+, 
\lb{2.21} \\
M_{\pm}(\bar z,x_0)=M_{\pm}( z,x_0)^*,\lb{2.22} \\
\rank (M_{\pm}(z,x_0))=m, \lb{2.23} \\
\lim_{\varepsilon\downarrow 0} M_{\pm}(\lambda
+i\varepsilon,x_0) \text{
exists for a.e.\
$\lambda\in\bbR$}. \lb{2.24}
\end{gather}
$\pm M_{\pm}(z,x_0)$ and $\mp M_{\pm}(z,x_0)^{-1}$ have 
isolated poles of
at most first order which
are real and have a negative definite residue.  \\
(ii)  $\pm M_{\pm}(z,x_0)$ admit the representations
\begin{align}
\pm M_{\pm}(z,x_0)&=F_\pm(x_0)+\int_\bbR 
d\Omega_\pm(\lambda,x_0) \,
\big(\f{1}{\lambda-z}-\f{\lambda}{1+\lambda^2}\big) 
\lb{2.25}  \\
&=\exp\bigg(C_\pm(x_0)+\int_\bbR d\lambda \, 
\Xi_{\pm} (\lambda,x_0)
\big(\f{1}{\lambda-z}-\f{\lambda}{1+\lambda^2}\big)  
 \bigg), \lb{2.26}
\end{align}
where
\begin{align}
F_\pm(x_0)&=F_\pm(x_0)^*, \quad \int_\bbR \f{\Vert
d\Omega_\pm(\lambda,x_0)\Vert}{1+\lambda^2}<\infty, 
\lb{2.27} \\
C_\pm(x_0)&=C_\pm(x_0)^*, 
\quad 0\le\Xi_\pm(\dott,x_0)\le I_m \, \rm{  a.e.}
\lb{2.28}
\end{align}
Moreover,
\begin{align}
\Omega_\pm((\lambda,\mu],x_0)&=\lim_{\delta\downarrow
0}\lim_{\varepsilon\downarrow 0}\f1\pi
\int_{\lambda+\delta}^{\mu+\delta} d\nu \, \Im(\pm
M_\pm(\nu+i\varepsilon,x_0)), \lb{2.29} \\
\Xi_\pm(\lambda,x_0)&=
\lim_{\varepsilon\downarrow 0}\f1\pi\Im(\ln(\pm
M_\pm(\lambda+i\varepsilon,x_0)) \text{ for a.e.\ 
$\lambda\in\bbR$}.\lb{2.30}
\end{align}
(iii)  Define the $n\times m$ matrices
\begin{align}
\Psi_\pm(z,x,x_0)&=\begin{pmatrix}\psi_{\pm,1}(z,x,x_0)\\
\psi_{\pm,2}(z,x,x_0)  \end{pmatrix} \no \\
&=\begin{pmatrix}\theta_1(z,x,x_0) & \phi_1(z,x,x_0)\\
\theta_2(z,x,x_0)& \phi_2(z,x,x_0)\end{pmatrix} 
\begin{pmatrix} I_m \\
M_\pm(z,x_0) \end{pmatrix}.
\lb{2.31}
\end{align}
Then the $m$ columns of $\Psi_\pm(z,x,x_0)$ form a basis for
$N(z,\pm\infty)$ and
\begin{equation}
\Im(M_\pm(z,x_0))=\Im(z) \int_{x_0}^{\pm\infty}dx\,
\Psi_\pm(z,x,x_0)^* A
\Psi_\pm(z,x,x_0) .\lb{2.32}
\end{equation}
\end{theorem}

In order to describe the Green's matrix associated with 
\eqref{2.4} on
$\bbR$, we assume the hypotheses of  Theorem~\ref{thm2.3} 
and introduce
\begin{align}
K(z,x,x^\prime)=\Psi_\mp(z,x,x_0)(M_-(z,x_0) &
-M_+(z,x_0))^{-1}\Psi_\pm(\bar z,x^\prime,x_0)^*, \no \\
& \hspace*{1.5cm} x\lessgtr x^\prime,\, 
z\in\bbC\backslash\bbR 
\lb{2.33}
\end{align}
and
\begin{align}
&M(z,x_0)=\f12(K(z,x_0,x_0+0)+K(z,x_0,x_0-0)) \lb{2.34} \\
&= \begin{pmatrix}N_-(z,x_0)^{-1}  
& \f12N_-(z,x_0)^{-1}N_+(z,x_0) \\
 \f12N_+(z,x_0)N_-(z,x_0)^{-1}
&M_\pm(z,x_0)N_-(z,x_0)^{-1}M_\mp(z,x_0)\end{pmatrix}, 
\no \\
& \hspace*{9cm} z\in\bbC\backslash\bbR, \lb{2.35}
\end{align}
with
\begin{equation}
N_{\pm}(z,x_0)=M_-(z,x_0) \pm M_+(z,x_0). \lb{2.35a}
\end{equation}
Next let $\phi\in L_A^2(\bbR)$ and consider
\begin{equation}
J\psi^\prime(z,x)=(zA+B(x))\psi(z,x)+A \phi(x), 
\quad z\in\bbC\backslash\bbR
\lb{2.36}
\end{equation}
for a.e.\ $x\in\bbR$.  Then \eqref{2.36} has a unique 
solution
$\psi(z,\dott)\in
L_A^2(\bbR)\cap\AC_{\loc}(\bbR)^n$ given by \cite{HS81}, 
\cite{HS83}
\begin{equation}
\psi(z,x)=\int_\bbR dx^\prime\, 
K(z,x,x^\prime) A\phi(x^\prime). \lb{2.37}
\end{equation}
Let $H$ be the matrix-valued Schr\"{o}dinger operator in
$L^2(\bbR)^m$,
\begin{align}
H&=-I_m \f{d^2}{dx^2}+Q, \lb{2.39} \\
\dom(H)&=\{g\in L^2(\bbR)^m \mid g,g^\prime
\in\AC_{\loc}(\bbR)^m,\,
 (-I_m g^{\prime\prime}+Qg)\in L^2(\bbR)^m\}. \no
\end{align}
In the following we associate the operator $H$ in 
$L^2(\bbR)^m$ with
the Hamiltonian system \eqref{2.4}. 
Hypothesis~\ref{hyp2.2} then renders
$H$ to be
self-adjoint in $ L^2(\bbR)^m$.
(Equivalently, the differential
expressions $-I_m d^2/dx^2+Q(x)$ is in the l.p.\ case at
$\pm\infty$.)

Denoting by $\rho(H)$, $\sigma(H)$, $\sigma_p(H)$, 
$\sigma_{\ess}(H)$,
$\sigma_{\ac}(H)$,
and $\sigma_{\singc}(H)$ the resolvent set, spectrum, 
point spectrum (i.e.,
the set of eigenvalues),
essential spectrum, absolutely and singularly continuous 
spectrum of $H$,
respectively, one can summarize the
connections between $M(z,x_0)$ and the various spectra of 
$H$ as follows.
\begin{theorem} [\cite{AD56}, \cite{Ca76}, \cite{GT97}, 
\cite{HS81},
\cite{HS82}, \cite{HS86}]
\lb{thm2.4} Assume Hypothesis~\ref{hyp2.2}, $z\in\bbC 
\backslash \bbR$, and
$x_0\in\bbR$.  Then \\
(i) $M(z,x_0)$ is a matrix-valued Herglotz function of 
rank $n$ with
representations
\begin{align}
M(z,x_0)&=F(x_0)+\int_\bbR d\Omega(\lambda,x_0)\,
\big(\f1{\lambda-z}-\f{\lambda}{1+\lambda^2}
\big)
\lb{2.42} \\
&=\exp\bigg(C(x_0)+\int_\bbR d\lambda\, 
\Upsilon (\lambda,x_0)
\big(\f1{\lambda-z}-\f{\lambda}{1+\lambda^2} \big) \bigg), 
\lb{2.43}
\end{align}
where
\begin{align}
F(x_0)&=F(x_0)^*, \quad \int_\bbR\f{\Vert 
d\Omega(\lambda,x_0)
\Vert}{1+\lambda^2}<\infty,
\lb{2.44}\\
C(x_0)&=C(x_0)^*, \quad 0\le \Upsilon (\dott,x_0)\le I_n 
\text{ a.e.} \lb{2.45}
\end{align}
Moreover,
\begin{align}
\Omega((\lambda,\mu],x_0)&=\lim_{\delta\downarrow
0}\lim_{\varepsilon\downarrow 0}\f1\pi
\int_{\lambda+\delta}^{\mu+\delta} d\nu \, \Im(
M(\nu+i\varepsilon,x_0)), \lb{2.46} \\
\Upsilon(\lambda,x_0)&=\lim_{\varepsilon\downarrow
0}\f1\pi\Im(\ln(M(\lambda+i\varepsilon,x_0)))
\text{ for a.e.\ $\lambda\in\bbR$}.\lb{2.47}
\end{align}
(ii) $z\in\rho(H)$ if and only if $M(z,x_0)$ is holomorphic 
near $z$.  In
this case,
\begin{equation}
((H-z)^{-1}f)(x)= \int_\bbR dx^\prime\, 
K_{1,1}(z,x,x^\prime)f(x^\prime),
\quad
z\in\bbC\backslash\sigma(H), \, f\in L^2(\bbR)^m. \lb{2.48} 
\end{equation}
{\rm (}Here $K_{1,1}(z,x,x^\prime)$ denotes the left upper 
$m\times m$ submatrix
of $K(z,x,x^\prime)$,
i.e., we write $K(z,x,x^\prime)=
(K_{j,k}(z,x,x^\prime))_{j,k=1}^2$.{\rm )}
Moreoever,
$\lambda_0\in\rho(H)\cap\bbR$ if and only if there is 
an $\varepsilon>0$
such that
\begin{equation}
\Omega(\lambda_0+\varepsilon,x_0)-\Omega(\lambda_0
-\varepsilon,x_0)=0.
\end{equation}
(iii) For all $\lambda\in\bbR$,
\begin{equation}
\lim_{\varepsilon\downarrow 0} \varepsilon\, 
\Im(M(\lambda+i\varepsilon,x_0))=
\Omega(\lambda+0,x_0)-\Omega(\lambda-0,x_0) \ge
0. \lb{2.51}
\end{equation}
(iv) $\lambda_0\in\sigma_p(H)$ if and only if 
\begin{equation}
\lim_{\varepsilon\downarrow 0} \varepsilon\, 
\Im(M(\lambda_0+i\varepsilon,x_0))=
\Omega(\lambda_0+0,x_0)-\Omega(\lambda_0-0,x_0) \neq
0. 
\lb{2.50}
\end{equation}
(v)
\begin{align}
\sigma(H)&=\supp(d\Omega(\dott,x_0)), \lb{2.52} \\
\sigma_{\ac}(H)&=\supp(d\Omega_{\ac}(\dott,x_0)), 
\lb{2.53a} \\
\sigma_{\singc}(H)&=\supp(d\Omega_{\singc}(\dott,x_0)), 
\lb{2.53b} \\
\overline{\sigma_p(H)}&=\supp(d\Omega_{\pp}(\dott,x_0)). 
\lb{2.54}
\end{align}
\end{theorem}
\noindent Here $\supp(\dott)$ denotes the topological 
(i.e., smallest closed)
support and
$d\Omega=d\Omega_{\pp}+d\Omega_{\singc}+d\Omega_{\ac}$ 
represents the
Lebesgue decomposition of
$d\Omega$ into its pure point (pp), singularly 
continuous (sc) and
absolutely continuous (ac) parts.

\section{Trace Formulas} \lb{s3}
In this section we derive trace formulas for the 
matrix-valued
Schr\"{o}dinger 
systems studied in Section~\ref{s2}. Throughout this 
section we will assume
the limit point case at
$\pm\infty$ and hence adopt Hypothesis~\ref{hyp2.2}.

In the Schr\"odinger case at hand, the inhomogeneous 
term in \eqref{2.36} is of the type 
$A\phi=(\phi_1,0)^t$ and $\psi=(\psi_1,\psi_2)^t=
(\psi_1,\psi_1^\prime)^t.$ Hence
$B(x)=\left(\begin{smallmatrix}-Q(x)& 0
\\ 0 & I_m
\end{smallmatrix}\right)$ is of a very
special nature and only the $m\times m$ submatrix 
$B_{1,1}=-Q(x)$ contains
information on $Q(x).$ Thus, 
we will focus on the $m\times m$ submatrix
$M_{1,1}(z,x_0)$  of $M(z,x_0)$ in
\eqref{2.35} and \eqref{2.42}--\eqref{2.45}.  By 
\eqref{2.33} and
\eqref{2.48} one infers that the
Green's matrix $G(z,x,x^\prime)$ of $H$ is given by
\begin{align}
G(z,x,x^\prime)&=K_{1,1}(z,x,x^\prime) \no \\
&=\psi_{\mp,1}(z,x,x_0)(M_-(z,x_0)-M_+(z,x_0))^{-1}
\psi_{\pm,1}(\bar z,x^\prime,x_0)^*, \no \\
& \hspace*{6cm} x\lesseqqgtr x^\prime, \, 
z\in\bbC\backslash\bbR, \lb{3.1}
\end{align}
where the $m\times m$ matrices $\psi_{\pm,j}(z,x,x_0)$ 
are defined in
\eqref{2.31}, that is,
\begin{equation}
\psi_{\pm,j}(z,x,x_0)=\theta_j(z,x,x_0)
+\phi_j(z,x,x_0)M_\pm(z,x_0), \quad
j=1,2. \lb{3.2}
\end{equation}
Moreover, since (block) diagonal elements of 
matrix-valued Herglotz
functions are (lower-dimensional) matrix-valued Herglotz 
functions (see, e.g., \cite{GT97}),
$M_{1,1}(z,x_0)$ $=G(z,x_0,x_0)$
is an $m\times m$ matrix-valued Herglotz function 
satisfying
\begin{align}
G(z,x_0,x_0)&=M_{1,1}(z,x_0)=(M_-(z,x_0)
-M_+(z,x_0))^{-1} \lb{3.3} \\
&= F_{1,1}(x_0)+\int_\bbR d\Omega_{1,1}(\lambda,x_0)\,
\big(\f1{\lambda-z}-\f{\lambda}{1+\lambda^2}
\big)  \lb{3.4}  \\
&=\exp\bigg(E(x_0)+\int_\bbR
d\lambda\, \Xi(\lambda,x_0)\,\big(\f1{\lambda-z}
-\f{\lambda}{1+\lambda^2}
\big)   \bigg), \lb{3.5}
\end{align}
where
\begin{align}
F_{1,1}(x_0)&=F_{1,1}(x_0)^*, \quad \int_\bbR\f{\Vert
d\Omega_{1,1}(\lambda,x_0)
\Vert}{1+\lambda^2}<\infty, \lb{3.6}  \\
E(x_0)&=E(x_0)^*, \quad 0\le\Xi(\lambda,x_0)\le I_m \,
 \rm{ a.e.,}  \lb{3.7}
\end{align}
and
\begin{align}
\Omega_{1,1}((\lambda,\mu],x_0)&=\lim_{\delta\downarrow 0}
\lim_{\varepsilon\downarrow
0}\f1\pi \int_{\lambda+\delta}^{\mu+\delta}d\nu\, 
\Im(G(\nu+i\varepsilon,x_0)),\lb{3.8} \\
\Xi(\lambda,x_0)&=\lim_{\varepsilon\downarrow 0}\f1\pi
\Im(\ln(G(\lambda+i\varepsilon,x_0))) \text{ for a.e.
$\lambda\in\bbR$}. \lb{3.9}
\end{align}

Next we discuss the asymptotic expansion of $G(z,x,x)$ as
$\abs{z}\to\infty$.  It will be
convenient to start with $M_\pm(z,x)$.  In the following 
we denote
$M(z)=O(\abs{z}^\alpha)$ or
$o(\abs{z}^\alpha)$ as $\abs{z}\to\infty$, whenever $\Vert
M(z)\Vert=O(\abs{z}^\alpha)$ or
$o(\abs{z}^\alpha)$ as $\abs{z}\to\infty$ for an 
appropriate matrix norm.

\begin{theorem} \cite[Theorem~4.7]{CG99} \lb{t3.1}
Suppose Hypothesis~\ref{hyp2.2}. In addition, assume  
$Q^{(N)}\in
L_{\loc}^1(\bbR)^{m\times m}$ for
some $N\in\bbN_0$ and let $C_\varepsilon\subset\bbC_+$ be 
the sector 
along the positive imaginary axis with vertex at zero 
and opening angle $\varepsilon$ with 
$0<\varepsilon<\pi/2$.  Then
$M_\pm(z,x)$ has the asymptotic
expansion as $\abs{z}\to\infty$ in $C_\varepsilon$ of 
the form
$(\Im(z^{1/2})\ge0$, $z\in\bbC)$
\begin{equation}
M_\pm(z,x)\underset{\substack{\abs{z}\to\infty\\ 
z\in C_\varepsilon}}{=}
\begin{cases}
\pm i I_m z^{1/2}+o(|z|^{1/2}) & \text{for $N=0$}, \\
\pm i I_m z^{1/2}+\sum_{k=1}^N m_{\pm,k}(x)z^{-k/2}
+o(|z|^{-N/2}) &
\text{for $N\in\bbN$}. \lb{3.10}
\end{cases}
\end{equation}
The expansion \eqref{3.10} is uniform with respect to 
$\arg\,(z)$ for $|z|
\to \infty$ in
$C_\varepsilon$ and uniform in $x\in\bbR$ as long as
$x$ varies in compact intervals.  The expansion 
coefficients $m_{\pm,k}(x)$
can be recursively
computed from
\begin{align}
m_{\pm,1}(x)&=\pm\f1{2i} Q(x), \quad m_{\pm,2}(x)=
 \f1{4} Q^\prime(x), \no \\
m_{\pm,k+1}(x)&=\pm\f{i}2\bigg(m_{\pm,k}^\prime(x)+
\sum_{\ell=1}^{k-1}m_{\pm,\ell}(x)m_{\pm,k-\ell}(x) 
\bigg), \quad k\ge 2.
\lb{3.11}
\end{align}
\end{theorem}

We briefly sketch a derivation of the recursion 
\eqref{3.11}. Let
\begin{equation}
\hat \Psi(z,x,x_0)=\begin{pmatrix} 
\psi_{-,1}(z,x,x_0)&\psi_{+,1}(z,x,x_0)\\
\psi_{-,2}(z,x,x_0)& \psi_{+,2}(z,x,x_0)\end{pmatrix}, 
\quad
z\in\bbC\backslash\bbR \lb{3.12}
\end{equation}
be the fundamental system of solutions of \eqref{2.14} 
as defined in
\eqref{2.31} and
observe that
\begin{equation}
\psi_{+,2}(z,x,x_0)=\psi_{+,1}^\prime(z,x,x_0) \lb{3.13}
\end{equation}
in the Schr\"{o}dinger operator case.  Hence any 
nonnormalized solutions
$\ti \psi_{\pm,1}(z,\dott)\in 
L^2((\pm\infty,c))^{m\times m}$, 
$c\in\bbR$ of
$-\psi_1^{\prime\prime}+Q\psi_1=z\psi_1$ for 
$z\in\bbC\backslash\bbR$ are of the
type
\begin{equation}
\ti \psi_{\pm,1}(z,x)=\psi_{\pm,1}(z,x,x_0)C_\pm \lb{3.14}
\end{equation}
for some nonsingular $m\times m$ matrices $C_\pm$.  Thus,
\begin{equation}
\ti \psi_{\pm,1}(z,x_0)=C_\pm, \quad \ti
\psi_{\pm,1}^\prime(z,x_0)=M_\pm(z,x_0)C_\pm \lb{3.15}
\end{equation}
by \eqref{2.17} and \eqref{2.31}.  In particular,
\begin{equation}
M_\pm(z,x_0)=\ti \psi_{\pm,1}^\prime(z,x_0) 
\ti \psi_{\pm,1}(z,x_0)^{-1}, \quad
z\in\bbC\backslash\bbR
\lb{3.16}
\end{equation}
is independent of the normalization chosen for 
$\ti \psi_{\pm,1}(z,x_0)$.
Varying the reference
point
$x_0\in\bbR$ then yields the standard Riccati-type 
equation,
\begin{equation}
M_\pm^\prime(z,x)+M_\pm(z,x)^2=Q(x)-z I_m. \lb{3.17}
\end{equation}
Existence of the asymptotic expansion \eqref{3.10} under 
the conditions
imposed on $Q(x)$ is a highly 
nontrivial matter and proved separately in \cite{CG99}. 
The recursion
relation for the
coefficients $m_{\pm,k}(x)$ in \eqref{3.11} then follows 
by inserting
\eqref{3.10} into \eqref{3.17}.

Since $G(z,x,x)=M_{1,1}(z,x)$, Theorem~\ref{t3.1} and 
\eqref{3.3} then yield an analogous
asymptotic expansion for the diagonal Green's matrix 
$G(z,x,x)$ of $H$.  In
fact, one obtains the
following result.

\begin{theorem} [\cite{CG99}] \lb{t3.2} 
Assume the hypotheses in Theorem~\ref{t3.1}.  Then 
$G(z,x,x)$ has an
asymptotic expansion
in $C_\varepsilon$ of the form
$(\Im(z^{1/2})\ge 0$ for $z\in\bbC)$
\begin{equation}
G(z,x,x)\underset{\substack{\abs{z}\to\infty\\ z\in
C_\varepsilon}}{=}\f{i}2 \sum_{k=0}^N G_k(x)
z^{-k-1/2}+o(|z|^{-N-1/2}), \lb{3.18}
\end{equation}
where
\begin{equation}
G_0(x)=I_m, \quad G_1(x)=\f12 Q(x), \text{ etc.} 
\lb{3.19}
\end{equation}
The expansion \eqref{3.18} is uniform with respect to 
$\arg\,(z)$ for $|z|
\to \infty$ in
$C_\varepsilon$ and uniform in $x\in\bbR$ as long as
$x$ varies in compact intervals.
\end{theorem}
\begin{proof}
The existence of the asymptotic expansion \eqref{3.18} 
is clear from
\eqref{3.3} and \eqref{3.10}. The actual expansion 
coefficients then can be determined from \eqref{3.3} 
and \eqref{3.11}.
\end{proof}

The trace formula for $Q(x)$ is then derived as follows.

\begin{theorem}  \lb{thm3.3}
In addition to Hypothesis~\ref{hyp2.2} suppose that
$Q\in\AC_{\loc}(\bbR)^{m\times m}$ and
$E_0=\inf(\sigma(H))>-\infty$.  Then
\begin{equation}
Q(x)=E_0 I_m+\lim_{z\to i\infty}\int_{E_0}^\infty 
d\lambda\,
z^2(\lambda-z)^{-2}(I_m-2\Xi(\lambda,x)), 
\quad x\in\bbR. \lb{3.21}
\end{equation}
\end{theorem}
\begin{proof}
By \eqref{3.5} and \eqref{3.9} one infers
\begin{equation}
\f{d}{dz}\ln(G(z,x,x))=\int_{E_0}^\infty d\lambda\, 
(\lambda-z)^{-2}
\Xi(\lambda,x).  \lb{3.22}
\end{equation}
By \eqref{3.18}, \eqref{3.19}, and the uniformity of 
the asymptotic expansion
\eqref{3.18} with respect to $\arg\,(z)$ as 
$|z| \to \infty$ in
$C_{\varepsilon}$, which
permits its differentiation in $z$, one derives
\begin{equation}
-\f{d}{dz}\ln(G(z,x,x))\underset{z\to i\infty}{=}
\f{1}{2}I_m z^{-1}+\f12
Q(x)z^{-2}+o(|z|^{-2}).
\lb{3.23}
\end{equation}
Thus,
\begin{align}
-\f{d}{dz}\ln(G(z,x,x))&=\f12 I_m(z-E_0)^{-1}
+\f12 \int_{E_0}^\infty d\lambda\,
(\lambda-z)^{-2}(I_m-2\Xi(\lambda,x)) \no \\
&\underset{z\to i\infty}{=}\f12 I_m z^{-1}
+ \f12 Q(x)z^{-2}+o(|z|^{-2})
\lb{3.24}
\end{align}
proves \eqref{3.21}.
\end{proof}

In the scalar case $m=1$, Theorem~\ref{thm3.3} was 
first derived in \cite{GS96}.
Subsequent extensions in the case $m=1$ and their 
applications to KdV invariants
 appeared in \cite{Ge95}, \cite{GH97}, \cite{GHSZ95}, 
\cite{Ry99}, \cite{Ry99a}. 
For an abstract approach to trace formulas based on 
perturbation theory and the theory of self-adjoint 
extensions of symmetric operators we refer to \cite{GM99}. 
The case of matrix-valued
Schr\"odinger operators was briefly sketched in 
\cite{GH97}. A different kind of trace
formula, based on scattering-theoretic concepts for 
short-range matrix-valued potentials,
appeared in \cite{MO82}. This reference also contains 
a variety of applications to
matrix-valued completely integrable evolution equations. 
Moreover, a trace formula for
matrix-valued Schr\"odinger operators $H$ on a finite 
interval with Dirichlet boundary
conditions at the endpoints was briefly discussed 
in \cite{Pa95}.

\begin{remark}\lb{rem3.4}
Assuming the hypotheses of Theorem~\ref{t3.1}, on infers
\begin{align}
-\f{d}{dz}\ln(G(z,x,x))&\underset{z\to i\infty}{=}
\sum_{k=0}^N R_k(x)
z^{-k-1}+o(|z|^{-N-1}), \no \\
R_0(x)&=\f12 I_m, \quad R_1(x)=\f12 Q(x), 
\text{ etc.}, \lb{3.25}
\end{align}
and derives in a similar fashion the higher-order trace 
formulas (see
\cite{GHSZ95} for the special scalar case $m=1$)
\begin{multline}
R_k(x)=\f12 E_0^k+k\,\lim_{z\to i\infty} 
\int_{E_0}^\infty d\lambda\,
z^{k+1}(\lambda-z)^{-k-1}(-\lambda)^{k-1}
(\f12 I_m-\Xi(\lambda,x)), \\
 k=1,\dots,N, \, x\in\bbR.\lb{3.26}
\end{multline}
\end{remark}

\section{Borg-Type Uniqueness Theorems} \lb{s4}
In 1946 Borg \cite{Bo46} proved, among a variety of other 
inverse spectral
theorems, the
following result.
\begin{theorem}[\cite{Bo46}] \lb{thm4.1}
Let $q\in L^2_{\loc} (\bbR)$ be real-valued and periodic. 
Let
$h=-\f{d^2}{dx^2}+q$ be the associated self-adjoint 
Schr\"odinger operator in $L^2(\bbR)$ $($cf.~\eqref{2.39} 
for $m=1)$ and suppose that
$\sigma(h)=[e_0,\infty)$ for some $e_0\in\bbR$.  Then
\begin{equation}
q(x)=e_0 \text{ for a.e. $x\in\bbR$}. \lb{4.1}
\end{equation}
\end{theorem}
Traditionally, uniqueness results such as 
Theorem~\ref{thm4.1} are called
Borg-type theorems.
(However, this terminology is not uniquely adopted 
and hence 
a bit
unfortunate. Indeed, inverse
spectral results on finite intervals recovering the 
potential
coefficient(s) from several spectra,
were also pioneered by Borg in his celebrated paper 
\cite{Bo46}, and hence
are also coined
Borg-type theorems in the literature, see, e.g., 
\cite[Sect.~6]{Ma94}.) The
purpose of this
section is to develop a new strategy of proof for such 
theorems based on
trace formulas and prove
extensions to matrix-valued situations.  In order to 
explain our strategy,
we provide a quick
proof of Theorem~\ref{thm4.1} (assuming
$q\in\AC_{\loc}(\bbR)$ instead of
$q\in L^2_{\loc} (\bbR)$).

\vspace*{3mm}
\noindent {\it Proof of Theorem~\ref{thm4.1}}
Suppose $q\in\AC_{\loc}(\bbR)$ is real-valued and periodic.  
Then by
standard Floquet theory, 
$g(z,x,x)$, the diagonal Green's function of $h$, is 
well-known to be
purely imaginary on
$\sigma(h)^o$ ($A^o$ the open interior of a set 
$A\subseteq\bbR$).  Hence, introducing 
$\xi(\lambda,x)=\pi^{-1}\lim_{\varepsilon\downarrow 0} 
\Im(\ln(g(\lambda+i\varepsilon,x,x))$ 
for~a.e.~$\lambda\in\bbR,$ one infers 
for all $x\in\bbR$,
\begin{equation}
\xi(\lambda,x)=\f12, \quad \lambda\in\sigma(h)^o. \lb{4.2}
\end{equation}
Since $\sigma(h)=[e_0,\infty)$ by hypothesis, the trace 
formula
\eqref{3.21} for $q(x)$ and
\eqref{4.2} yield
\begin{equation}
q(x)=e_0+\lim_{z\to i\infty}\int_{e_0}^\infty 
d\lambda\, z^2
(\lambda-z)^{-2}(1-2\xi(\lambda,x))=e_0, 
\quad x\in\bbR. \lb{4.3}
\end{equation}
\hspace*{12.2cm} $\square$

A closer examination of the proof shows that periodicity 
of $q(x)$ is not
the point for the
uniqueness result \eqref{4.1}.  The key ingredient (besides
$\sigma(h)=[e_0,\infty)$ and $q$ real-valued) is clearly
the fact
\eqref{4.2}, that is, for all $x\in\bbR$,
\begin{equation}
\xi(\lambda,x)=1/2 \text{ for a.e. } 
\lambda\in\sigma_{\ess}(h)
\lb{4.3a}
\end{equation}
($\sigma_{\ess}(\dott)$ the essential spectrum).

Real-valued periodic potentials are known to satisfy 
\eqref{4.3a} but so
are certain classes of real-valued
quasi-periodic and almost-periodic potentials 
$q(x)$ (see, e.g.,
\cite{BE95}, 
\cite{Cr89}, \cite{DS83}, \cite{Ko84}, \cite{Ko87a}, 
\cite{Ko87b},
\cite{KK88}, \cite{KS88}, \cite{SY95}).
In particular, the class of real-valued algebro-geometric 
finite-gap potentials
$q(x)$ (a subclass of the set of real-valued quasi-periodic 
potentials) is
a prime example
satisfying \eqref{4.3a} without necessarily being periodic. 
Traditionally,
potentials $q(x)$
satisfying \eqref{4.3a} are called \textit{reflectionless} 
(see
\cite{BE95}, 
\cite{Cr89}, \cite{DS83}, \cite{KK88}).

\begin{remark} \lb{r4.1a}
We note that real-valuedness of $q$ is an essential 
assumption in Theorem~\ref{thm4.1}. Indeed, $q(x)=\exp(ix),$ 
$x\in\bbR,$ is well-known to lead to the half-line spectrum 
$\sigma (h)=[0,\infty),$ with $h=-\f{d^2}{dx^2}+q$ in 
$L^2(\bbR)$ defined on the standard Sobolov space 
$H^{2,2}(\bbR).$ A detailed treatment of a class of examples 
of this type can be found in \cite{Ga80}, \cite{Ga80a}, 
\cite{GU83}, \cite{PT88}, \cite{PT91}. Moreover, the example 
of complete exponential localization of the spectrum of a 
discrete Schr\"odinger operator with a quasi-periodic 
real-valued potential 
having two basic fequencies and no gaps in its spectrum 
illustrates the importance of the reflectionless property of 
$q$ in Theorem~\ref{thm4.1}.
\end{remark}

Taking the quick proof of Theorem~\ref{thm4.1} as the 
point of departure
for our extension of
Borg-type results to matrix-valued Schr\"{o}dinger 
operators, we
now use the reflectionless
situation described in \eqref{4.3a} as the model for 
the subsequent
definition.

\begin{definition}\lb{def4.2}
Assume Hypothesis~\ref{hyp2.2} and 
$Q\in\AC_{\loc}(\bbR)^{m\times m}$. \\
Then the matrix-valued potential
$Q(x)$ is called {\it reflectionless} if for all 
$x\in\bbR$,
\begin{equation}
\Xi(\lambda,x)=\f12 I_m 
\text{  for a.e.\ $\lambda\in\sigma_{\ess}(H)$}.
\lb{4.5}
\end{equation}
\end{definition}
Since hardly any confusion can arise, we will also 
call $H$ 
reflectionless if \eqref{4.5} is satisfied. 

Given Definition~\ref{def4.2}, we turn to a Borg-type 
uniqueness theorem and formulate the analog of
Theorem~\ref{thm4.1} for
(reflectionless) matrix-valued Schr\"{o}dinger operators.
\begin{theorem}\lb{thm4.6}
Assume Hypothesis~\ref{hyp2.2} and 
$Q\in\AC_{\loc}(\bbR)^{m\times m}$. 
Suppose that
$Q(x)$ is reflectionless and $\sigma(H)=[E_0,\infty)$. 
Then
\begin{equation}
Q(x)=E_0 I_m \text{ for all $x\in\bbR$}.  \lb{4.35}
\end{equation}
\end{theorem}
\begin{proof}
By hypothesis, $\Xi(\lambda,x)=(1/2)I_m$ for a.e.\
$\lambda\in[E_0,\infty)$ and all
$x\in\bbR$.  Thus the trace formula \eqref{3.21} 
yields \eqref{4.35}.  
\end{proof}

In the remainder of the section we will show that the 
case of periodic
$Q(x)$ is covered by Theorem~\ref{thm4.6} under appropriate 
uniform multiplicity assumptions on
$\sigma(H)$.
Among other results this then recovers a recent theorem 
by D\'{e}pres \cite{De95} for
matrix-valued periodic Schr\"{o}dinger operators.

In order to discuss Floquet theory for $H$ we adopt the 
following
assumptions for the remainder of this section.

\begin{hypothesis} \lb{hyp4.7}
In addition to Hypothesis~\ref{hyp2.1} suppose that 
$Q\in\AC_{\loc}(\bbR)^{m\times m}$ is periodic, that
is, there is an $\omega>0$ such that $Q(x+\omega)=Q(x)$ 
for all $x\in\bbR$.
\end{hypothesis}

Since by Hypothesis~\ref{hyp4.7}, 
$Q\in L^{\infty}(\bbR)^{m \times m}$,
the corresponding periodic
Hamiltonian system
\eqref{2.4} is in the l.p. case at $\pm \infty$.

We briefly review a few basic facts from Floquet theory 
for Hamiltonian
systems of the type
$J\psi^\prime(z,x)=(zA+B(x))\psi(z,x)$ with $B(x)$ 
periodic of period
$\omega>0$. For a detailed
treatment of Floquet theory, relevant in our context, 
we refer to 
\cite{De80}, \cite[pp.~1486--1498]{DS88}, \cite{GL87},
\cite{Ha80}, \cite{Kh77}, 
\cite{Kr83}, \cite{Kr83a}, \cite{MV81},
\cite{MV89}, \cite{Ro63},
\cite{Ve83}, \cite{We87}, \cite{Ya92}, \cite{YS75a}, 
\cite{YS75b}, and the
literature therein. Recalling the
notation introduced in \eqref{2.16} and \eqref{2.17} one 
considers the
monodromy matrix
\begin{equation}
\Phi(z,x_0)=\begin{pmatrix} \theta_1(z,x_0+\omega,x_0) &
\phi_1(z,x_0+\omega,x_0)\\
 \theta_2(z,x_0+\omega,x_0) &  
\phi_2(z,x_0+\omega,x_0)\end{pmatrix}, \quad
z\in\bbC. \lb{4.37}
\end{equation}
Denoting its eigenvalues by $\rho_j(z)$, that is,
\begin{equation}
\sigma(\Phi(z,x_0))=\{\rho_j(z)\}_{j=1,\dots,n}, \lb{4.38}
\end{equation}
it is a well-known fact that $\sigma(\Phi(z,x_0))$, 
unlike $\Phi(z,x_0)$,
is independent of the
chosen reference point $x_0\in\bbR$. Moreover,
\begin{equation}
\det(\Phi(z,x_0))=1, \quad z \in \bbC. \lb{4.38a}
\end{equation}
One then obtains the following characterization of the 
spectrum of $H$,
\begin{equation}
\sigma(H)=\{\lambda\in\bbR \mid \abs{\rho_j(\lambda)}=1
\text{ for some $j\in\{1,2,\dots,n \}$} \}. \lb{4.39}
\end{equation}
In particular,
\begin{equation}
\abs{\rho_j(z)}\neq 1 \text{ for all 
$z\in\bbC\backslash\sigma(H)$}.
\lb{4.40}
\end{equation}
Let $\Psi(z,x,x_0)$ denote the normalized fundamental 
system \eqref{2.16}
of the
Hamiltonian system \eqref{2.14} with 
$\Psi(z,x_0,x_0)=I_n$, then
\begin{equation}
\Psi(z,x+\omega,x_0)=\Psi(z,x,x_0)\Phi(z,x_0), 
\quad z\in\bbC \lb{4.40a}
\end{equation}
by periodicity of $Q$. Since by hypothesis $H$ is 
self-adjoint, there
exists a fundamental
system $\hat \Psi(z,x,x_0)$ of \eqref{2.14} of the 
following type (cf.~\eqref{2.31}),
\begin{equation}
\hat \Psi(z,x,x_0)=
\begin{pmatrix} \psi_{-,1}(z,x,x_0) 
& \psi_{+,1}(z,x,x_0) \\
\psi_{-,2}(z,x,x_0)& \psi_{+,2}(z,x,x_0)\end{pmatrix}, 
\quad
z\in\bbC\backslash\sigma(H), \lb{4.41}
\end{equation}
with $\psi_{\pm,1}(z,x_0,x_0)=I_m$, where
\begin{equation}
\psi_{\pm,1}(z,\dott,x_0)\in L^2((R,\pm\infty))^{m\times m} 
\text{ for all
$R\in\bbR$},
\, \, z\in\bbC \backslash \sigma(H). \lb{4.42}
\end{equation}
Thus,
\begin{equation}
\hat \Psi(z,x+\omega,x_0)= \hat \Psi(z,x,x_0) 
\hat \Phi(z,x_0), \quad
z\in\bbC \backslash \sigma(H) \lb{4.43}
\end{equation}
and hence, as in the proof of Theorem~\ref{t3.1} (cf.\
\eqref{3.15}), one infers that $\hat \Phi(z,x_0)$ must 
be of the form
\begin{equation}
\hat \Phi(z,x_0)=\begin{pmatrix}\rho_-(z,x_0) & 0 \\ 0 &
\rho_+(z,x_0)\end{pmatrix} \lb{4.44}
\end{equation}
for nonsingular $m\times m$ matrices $\rho_{\pm}(z,x_0)$.  
Thus,
\begin{equation}
\psi_{\pm,j}(z,x+\omega,x_0)=
\psi_{\pm,j}(z,x,x_0)\rho_\pm(z,x_0), \quad
j=1,2. \lb{4.45}
\end{equation}
Next, noticing
\begin{equation}
\hat \Psi(z,x,x_0) = \Psi(z,x,x_0)C(z,x_0), 
\quad z\in\bbC \backslash
\sigma(H) \lb{4.45a}
\end{equation}
for some nonsingular $n \times n$ matrix $C(z,x_0)$ one 
infers
\begin{equation}
\hat \Phi(z,x_0) = C(z,x_0)^{-1} \Phi(z,x_0) C(z,x_0), \quad
z\in\bbC \backslash \sigma(H) \lb{4.45b}
\end{equation}
and hence
\begin{equation}
\sigma (\hat \Phi(z,x_0)) = \sigma (\Phi(z,x_0)), \quad  
z\in\bbC \backslash
\sigma(H). \lb{4.45c}
\end{equation}
We observe from \eqref{4.42}, \eqref{4.44}, \eqref{4.45}, 
and
\eqref{4.45c} that $\sigma(\Phi(z,x_0))$ can be 
partioned as
\begin{align}
\sigma(\Phi(z,x_0))=\{ & \rho_j(z) \}_{j=1,\dots,n}
=\sigma(\rho_-(z,x_0))\cup\sigma(\rho_+(z,x_0)), \lb{4.52}
\\ & \sigma(\rho_\pm(z,x_0))=
\{\rho_{\pm,j}(z) \}_{j=1,\dots,m}, \lb{4.52a}
\end{align}
where (cf. \eqref{4.38a} and \eqref{4.40})
\begin{equation}
0 \neq \abs{\rho_{\pm,j}(z)}\lessgtr 1 \text{ for
$z\in\bbC\backslash\sigma(H)$, $j=1,\dots, m$}.
\lb{4.53}
\end{equation}
Hence
\begin{equation}
\psi_{\pm,1}(z,x_0+\omega,x_0)=\rho_\pm(z,x_0)=
\theta_1(z,x_0+\omega,x_0)+
\phi_1(z,x_0+\omega,x_0)M_\pm(z,x_0) \lb{4.46}
\end{equation}
yields
\begin{align}
M_\pm(z,x_0)&=\psi_{\pm,2}(z,x_0,x_0) \no \\
&=\phi_1(z,x_0+\omega,x_0)^{-1}\big(\rho_\pm(z,x_0)-
\theta_1(z,x_0+\omega,x_0) 
\big), \lb{4.47} \\
&\hspace*{4.5cm} z\in\bbC\backslash(\sigma(H)
\cup\sigma(H_{x_0}^D)) \no
\end{align}
and
\begin{align}
\psi_{\pm,j}(z,x,x_0)=&\theta_j(z,x,x_0) \no \\
&
+\phi_j(z,x,x_0)\phi_1(z,x_0+
\omega,x_0)^{-1}(\rho_\pm(z,x_0)-
\theta_1(z,x_0+\omega,x_0)),  \no \\
& \hspace*{3.5cm} 
z\in\bbC\backslash(\sigma(H)\cup\sigma(H_{x_0}^D)), \,
j=1,2, \lb{4.48}
\end{align}
where
\begin{equation}
\sigma(H_{x_0}^D)=
\{z\in\bbC \mid \det(\phi_1(z,x_0+\omega,x_0))=0 \}.
\lb{4.49}
\end{equation}
One can show that $\sigma(H_{x_0}^D)$ is the spectrum 
of a self-adjoint
operator
$H_{x_0}^D$ (associated with a Dirichlet-type boundary 
condition
$\psi_1(z,x_0+\omega)=\psi_1(z,x_0)=0$, 
$\psi_1(z,\dott)\in \AC([x_0,x_0
+\omega])^m$) and
hence
$\sigma(H_{x_0}^D)\subset\bbR$.  Combining
\eqref{2.16}, \eqref{4.37}, \eqref{4.40a}, \eqref{4.45}, 
 and \eqref{4.48} yields
\begin{multline}
\hat \Phi(z,x_0)^2-\begin{pmatrix}
\theta_1+\phi_1\phi_2\phi_1^{-1} &
0 \\ 0 &
\theta_1+\phi_1\phi_2\phi_1^{-1}
\end{pmatrix}\hat \Phi(z,x_0) \\
+\begin{pmatrix}\phi_1\phi_2\phi_1^{-1}\theta_1
-\phi_1\theta_2 & 0 \\ 0 &
\phi_1\phi_2\phi_1^{-1}\theta_1
-\phi_1\theta_2\end{pmatrix}=0, \lb{4.50}
\end{multline}
where $\theta_j, \phi_j$ are evaluated at the point 
$(z,x_0+\omega,x_0)$,
that is,
\begin{equation}
\phi_j=\phi_j(z,x_0+\omega,x_0), \quad \theta_j
=\theta_j(z,x_0+\omega,x_0),
\quad j=1,2. \lb{4.50a}
\end{equation}
Equation \eqref{4.50} is  equivalent to
\begin{align}
& \rho_{\pm}(z,x_0)^2 -\big(\theta_1(z,x_0+\omega,x_0)+
\phi_1(z,x_0+\omega,x_0)\phi_2(z,x_0+\omega,x_0)\times 
\no \\
& \hspace*{2.05cm}\times\phi_1(z,x_0+\omega,x_0)^{-1}\big)
\rho_{\pm}(z,x_0) \no \\
&+\phi_1(z,x_0+\omega,x_0)
\phi_2(z,x_0+\omega,x_0)\phi_1(z,x_0+\omega,x_0)^{-1}
\theta_1(z,x_0+\omega,x_0) \lb{4.51} \\
& -\phi_1(z,x_0+\omega,x_0)\theta_2(z,x_0+\omega,x_0)=0, 
\quad
z\in\bbC\backslash(\sigma(H)\cup\sigma(H_{x_0}^D)). \no
\end{align}
In anticipation of \eqref{4.57}, which will be proven next,  
$\rho_{\pm}(z,x_0)$ can be extended to
$\rho_{\pm}(\lambda+i0,x_0)=\lim_{\varepsilon\downarrow
0}\rho_{\pm}(\lambda+i\varepsilon,x_0)$ by continuity 
for all
$\lambda\in\sigma(H)^o.$ Hence we will in the following 
extend the domain of validity of \eqref{4.44} and 
\eqref{4.50}--\eqref{4.51} to all $z\in\sigma(H)^o$ 
(agreeing to take normal limits to the real line in 
$\bbC_+$).

\begin{theorem} \lb{thm4.8}
Suppose Hypothesis~\ref{hyp4.7}.  If $H$ has uniform
spectral multiplicity $2m$, then for all $x\in\bbR$ 
and all $\lambda\in\sigma(H)^o,$
\begin{equation}
M_+(\lambda+i0,x)=M_-(\lambda+i0,x)^*=M_-(\lambda-i0,x). 
\lb{4.55}
\end{equation}
In particular, $M_-(z,x)$ is the analytic continuation of 
$M_+(z,x)$ (and vice versa) through $\sigma(H)^o$.
\end{theorem}
\begin{proof}
By general Floquet theory, $\sigma(H)$ consists of a 
countable union of
closed intervals on $\bbR$,
possibly separated by gaps in between. Moreover, $H$ is 
bounded from
below and it has no eigenvalues,
\begin{equation}
\sigma_p(H)=\emptyset, \quad \sigma(H)=\sigma_{\ess}(H)
=\sigma_c(H). \lb{4.56}
\end{equation}
In fact, one can show that $\sigma(H)$ is purely 
absolutely continuous,
$\sigma(H)=\sigma_{\ac}(H)$ (as is also clear from the 
analytic continuation 
of $M_\pm(\lambda +i0,x)$ through $\sigma(H)^o$ implied 
by \eqref{4.55}),
but we omit the details.  The assumptions of uniform 
(maximal) spectral multiplicity $n=2m$ of
$\sigma(H)$ guarantees the existence of $n$ eigenvalues 
$\rho_j(\lambda)$
of the monodromy matrix
$\Phi(\lambda,x_0)$ with $\abs{\rho_j(\lambda)}=1$ for 
$j=1,\dots,n$ for
$\lambda\in\sigma(H)^o$,
in particular, $\Phi(\lambda,x_0)$ is unitary and hence 
diagonalizable for
$\lambda\in
\sigma(H)^o$. 

Next, suppose that $\psi_1(\lambda,\dott,x_0),\dots,
\psi_n(\lambda,\dott,x_0)  \in
L^{\infty}(\bbR)^n$  for
$\lambda\in\sigma(H)^o$, with
$\psi_j(\lambda,x_0,x_0)=(\delta_{j,1},
\dots,\delta_{j,n})^t$ for
$j=1,\dots,n$ are $n$ linearly
independent normalized solutions of \eqref{2.4}. We 
claim that
\begin{equation}
\sigma_p(H_{x_0}^D)\cap\sigma(H)^o=\emptyset  \lb{4.57}
\end{equation}
since eigenfunctions of $H_{x_0}^D$ for 
$\lambda\in\sigma(H)^o$ would
necessarily be linear
combinations of $\psi_1(\lambda,x,x_0),\dots, 
\psi_n(\lambda,x,x_0)$.
However, none of them can lie in
$L^2(\bbR)^n$ since the fundamental matrix
$\Psi(\lambda,x,x_0)=(\psi_1(\lambda,x,x_0),\dots, 
\psi_n(\lambda,x,x_0))$,
$\lambda\in\sigma(H)^o$ satisfies
\begin{equation}
\Psi(\lambda,x+\omega, x_0)=
\Psi(\lambda,x,x_0)\Phi(\lambda,x_0), \quad
\lambda\in\sigma(H)^o \lb{4.58}
\end{equation}
(cf.\ \eqref{4.40a}) with $\Phi(\lambda,x_0)$ unitary. 
Since
$\sigma(H_{x_0}^D)\subset\bbR$, 
\eqref{4.57} implies the existence of 
$\phi_1(z,x_0+\omega,x_0)^{-1}$ for
$z\in\bbC\backslash \sigma(H_{x_0}^D)$. 

In the following we denote by
$\calD\subset\bbC$ the
discrete set of points (i.e., countable without finite 
limit points)  where
$\hat \Phi(z,x_0)$ is
not diagonalizable. 

By our hypothesis of uniform spectral 
multiplicity $n$ of $H$,
$\hat \Phi(z,x_0)$ is diagonalizable for all
$z\in\bbC\backslash \partial \sigma(H)$ ($\partial A$ 
denoting the boundary of
a  subset $A\subseteq \bbR$) and we conclude that 
$\calD\subseteq
\partial \sigma(H)$. By a well-known argument, see, 
for instance
\cite[Lemma XIII.7.63]{DS88}, 
one infers that all points in $\partial \sigma(H)$ are 
branch points for
the eigenvalues
$\rho_j(z)$ of $\Phi(z,x_0)$ and hence 
$\calD=\partial \sigma(H)$. Thus
we obtain upon diagonalizing
$\hat \Phi(z,x_0)$ in
\eqref{4.50} that
\begin{equation}
\ti\Phi(z,x_0)^2-E(z,x_0)\ti\Phi(z,x_0)+F(z,x_0)=0, 
\quad z\in\bbC\backslash
(\partial \sigma(H) \cup \sigma(H^D_{x_0})),
\lb{4.59}
\end{equation}
where
\begin{equation}
\ti\Phi(z,x_0)=\begin{pmatrix}
\rho_1(z) & 0 & \dots & 0 \\
0 & \rho_2(z) & \dots & 0 \\
\vdots & \vdots & \ddots &  \vdots \\
0& 0 & \dots & \rho_n(z)
\end{pmatrix} \lb{4.60}
\end{equation}
denotes the diagonalization of $\hat \Phi(z,x_0)$, 
and $E(z,x_0)$,
$F(z,x_0)$ are analytic in
$z\in\bbC\backslash\sigma(H)$ with a continuous 
extension to
$\sigma(H)^o$.  Thus each
$\rho_j(z)$  satisfies a quadratic equation and one 
introduces a canonical
set of cuts
$\{\calC_k\}_{k\in I}$ along $\sigma(H)$ ($I$ a 
finite or countably
infinite index set), joining the branch points, that is, 
all points in
$\partial \sigma(H)$ (as well as $+\infty$ and possibly 
$-\infty$ in case
$I$ is finite)
\begin{equation}
\sigma(H) = \bigcup_{k \in I}\calC_k. \lb{4.60a}
\end{equation}
In this manner, each
$\rho_j(z)$ becomes an  analytic function on a (fixed) 
two-sheeted Riemann
surface (glued together
crosswise along these cuts in a standard manner). In 
particular,
$\{\rho_1(z),\dots, \rho_n(z) \}$ can now be split 
into pairs
$\{\rho_{1,+}(z), \rho_{1,-}(z), \dots,$ $\rho_{m,+}(z), 
\rho_{m,-}(z) \}$
such
that $\rho_{k,-}(z)$ represents the analytic continuation 
of $\rho_{k,+}(z)$ 
(and vice versa), whenever $z$ crosses tranversally 
through one of the cuts.  

Next, pick a $\lambda_0\in\sigma(H)^o$, that is, 
$\lambda_0\in\calC_{k_0}^o$ for
some
$k_0\in I$ and pick a
$\rho_{j_0,+}(z)$ for $z$ in a sufficiently small 
neighborhood
$U_+(\lambda_0)\cap\bbC_+$ (or
$U_-(\lambda_0)\cap\bbC_-)$ of
$\lambda_0$. Suppose $\abs{\rho_{j_0,+}(z)}>1$ for 
$z$ along a path in
$U_+(\lambda_0)$ transversally
approaching $\calC_{k_0}^0$ and intersecting 
$\calC_{k_0}^o$ at
$\lambda_0$.  By analyticity of
$\rho_{j_0,+}(z)$, the analytic continuation $\ti
\rho_{j_0,+}(z)=\rho_{j_0,-}(z)$ of
$\rho_{j_0,+}(z)$ will satisfy 
$\abs{\ti \rho_{j_0,+}(z)}<1$ in an
appropriate neighborhood
$V_-(\lambda_0)\cap\bbC_-$ (or $V_+(\lambda_0)\cap\bbC_+$) 
of $\lambda_0$.
Hence we may identify
$\{\rho_{j,\pm}(z)\}_{j=1,\dots,m}$ with
$\{\rho_{\pm,j}(z)\}_{j=1,\dots,m}$  in
\eqref{4.52}--\eqref{4.53}.  Thus upon possibly reordering 
the eigenvalues
along the diagonal in
\eqref{4.60} we may write
\begin{equation}
\ti\Phi(z,x_0)=\begin{pmatrix}\ti\rho_-(z) & 0 \\ 0 &  
\ti \rho_+(z)
\end{pmatrix},
\lb{4.61}
\end{equation}
where
\begin{equation}
\ti\rho_\pm(z)=\begin{pmatrix}
\rho_{\pm,1}(z)& 0 & \dots & 0 \\
 0 &       \rho_{\pm,2}(z)& \dots & 0\\
 \vdots& \vdots &\ddots & \vdots  \\
 0& 0& \dots &\rho_{\pm,m}(z)\end{pmatrix}
\lb{4.62}
\end{equation}
such that $\ti\rho_-(z)$ is the analytic continuation of 
$\ti\rho_+(z)$
through the
interior of the cuts $\cup_{k\in I}\calC_k^o$ and hence 
through
$\sigma(H)^o$. Since the
similarity transformations connecting $\ti \rho_{\pm}(z)$ 
and
$\rho_{\pm}(z,x_0)$ can be
chosen as matrices whose column vectors are the 
eigenvectors $e_{\pm
,j}(z)$ of
$\rho_{\pm}(z,x_0)$, and $e_{\pm ,j}(z)$ have the same 
branching behavior
as the associated
eigenvalues $\rho_{\pm ,j}(z)$ (see, e.g.,
\cite[Sect.~6.1]{Ba85}),  one infers that
$\rho_-(z,x_0)$ is the analytic continuation of 
$\rho_+(z,x_0)$ through
$\sigma(H)^o$. By
\eqref{4.47},
$M_-(z,x_0)$ is the analytic continuations of
$M_+(z,x_0)$ through $\sigma(H)^o$. Consequently,
\begin{equation}
M_-(\lambda\mp i0,x_0)=M_-(\lambda\pm i0,x_0)^*=
M_+(\lambda\pm i0,x_0),
\quad \lambda\in \sigma(H)^o. \lb{4.63}
\end{equation}
Since $x_0\in\bbR$ was arbitrary we obtain \eqref{4.55}.
\end{proof}

The next result proves
necessary and sufficient conditions for $Q$
 to be reflectionless. It is modeled after Lemma~3.3 
in \cite{GKT96} in the context of Jacobi operators and we 
provide a formulation that anticipates extensions to the 
non-periodic case following \cite{SY95}, to be discussed 
elsewhere. 

\begin{theorem} \lb{thm4.3}
Assume Hypothesis~\ref{hyp4.7} and that $H$ has uniform 
spectral multiplicity $2m.$ Let 
$\Sigma\subseteq\sigma(H)^o.$ Then the
following conditions are equivalent and each of them 
holds. \\
(i) For all $x\in\bbR$ and all $\lambda\in\Sigma,$
\begin{equation}
\Xi(\lambda,x)=\f12 I_m.
\lb{4.7}
\end{equation}
(ii) For some $x_0\in\bbR$ and all $\lambda\in\Sigma,$
\begin{align}
G(\lambda+i0,x_0,x_0)&=
-G(\lambda+i0,x_0,x_0)^*, \lb{4.8} \\
G^\prime(\lambda+i0,x_0,x_0)&=
-G^\prime(\lambda+i0,x_0,x_0)^*.
\lb{4.9}
\end{align}
(Here
$G^\prime(\lambda+i0,x_0,x_0)=
\f{d}{dx}G(\lambda+i0,x,x)|_{x=x_0}$.)
\\
(iii)  For some $x_0\in\bbR$ and all $\lambda\in\Sigma,$
\begin{equation}
M_+(\lambda+i0,x_0)=M_-(\lambda+i0,x_0)^*. \lb{4.10}
\end{equation}
\end{theorem}
\begin{proof}
In the following, let $x,x_0\in\bbR$ and $\lambda\in\Sigma.$ 
By Theorem~\ref{thm4.8}, the normal limits
$M_\pm(\lambda+i0,x_0)$  and
$(M_-(\lambda+i0,x_0)-M_+(\lambda+i0,x_0))^{-1}$ 
exist for all $\lambda\in\Sigma$ and condition (iii) holds. 
By \eqref{3.9}, \eqref{4.7} is equivalent to
\begin{equation}
G(\lambda+i0,x,x)=-G(\lambda+i0,x,x)^*  \lb{4.12}
\end{equation}
and hence taking $x=x_0$ implies \eqref{4.8}. By 
\eqref{2.16}, \eqref{3.1}, and \eqref{3.2}, one infers
\begin{align}
G(z,x,x)&=\psi_{\mp,1}(z,x,x_0)G(z,x_0,x_0)
\psi_{\pm,1}(\bar z,x,x_0)^* \no \\
&=(\theta_1(z,x,x_0)+\phi_1(z,x,x_0)M_\mp(z,x_0))G(z,x_0,x_0)
\times \no \\
&\quad \, \times 
(\theta_1(z,x,x_0)+M_\pm(z,x_0)\phi_1(z,x,x_0)), 
\quad z\in\bbC_+, \lb{4.12a}
\end{align}
and 
\begin{align}
\f{d}{dx}G(z,x,x) 
&=(\theta_1^\prime (z,x,x_0)+\phi_1^\prime (z,x,x_0)
M_\mp(z,x_0))G(z,x_0,x_0) \times \no \\
&\quad \, \times 
(\theta_1 (z,x,x_0)+M_\pm(z,x_0)
\phi_1 (z,x,x_0)) \no \\
&+(\theta_1 (z,x,x_0)+\phi_1 (z,x,x_0)
M_\mp(z,x_0))G(z,x_0,x_0) \times \no \\
&\quad \, \times 
(\theta_1^\prime (z,x,x_0)+M_\pm(z,x_0)
\phi_1^\prime (z,x,x_0)), \quad z\in\bbC_+. \lb{4.12b}
\end{align}
Hence, $G(\lambda+i0,x,x)$ and 
$G^\prime (\lambda+i0,x,x)$ exist for all 
$\lambda\in\Sigma,$ and we
may  differentiate \eqref{4.12}
with respect to $x\in\bbR$ to obtain 
\begin{equation}
\f{d}{dx}G(\lambda+i0,x,x)=
-\f{d}{dx}G(\lambda+i0,x,x)^*.  \lb{4.12c}
\end{equation}
Taking $x=x_0$ in \eqref{4.12c} then implies \eqref{4.9} 
and hence we have shown that (i) implies (ii).  

Next we prove that (ii) implies (iii). We could immediately 
invoke \eqref{4.55}, but prefer to show a simple argument 
that permits extensions to non-periodic cases to be discussed 
elsewhere. By \eqref{3.3}, \eqref{4.12} is equivalent to
\begin{equation}
M_-(\lambda+i0,x)-M_+(\lambda+i0,x)=-(M_-(\lambda+i0,x)
-M_+(\lambda+i0,x))^*
\lb{4.13}
\end{equation}
and hence to
\begin{equation}
\Re(M_+(\lambda+i0,x))=\Re(M_-(\lambda+i0,x)). \lb{4.14}
\end{equation} 
Thus, \eqref{4.8} is equivalent to
\begin{equation}
\Re(M_+(\lambda+i0,x_0))=\Re(M_-(\lambda+i0,x_0)) 
\text{ for all\
$\lambda\in\Sigma$}. \lb{4.15}
\end{equation}
In order to exploit \eqref{4.9} one computes from 
\eqref{3.3} and \eqref{4.12b},
\begin{align}
G^\prime(z,x_0,x_0) = \, &M_\mp(z,x_0)(M_-(z,x_0)-
M_+(z,x_0))^{-1} \no\\
&+(M_-(z,x_0)-M_+(z,x_0))^{-1}M_\pm(z,x_0). \lb{4.16}
\end{align}
Consequently, \eqref{4.9} is equivalent to
\begin{multline}
(M_\mp(\lambda+i0,x_0)-M_\pm(\lambda+i0,x_0)^*)
(M_-(\lambda+i0,x_0)-M_+(\lambda+i0,x_0))^{-1} \\
=(M_-(\lambda+i0, x_0)-
M_+(\lambda+i0,x_0))^{-1}(M_\mp(\lambda+i0,x_0)^*-
M_\pm(\lambda+i0,x_0)). \lb{4.17}
\end{multline}
A simple manipulation in \eqref{4.17} (adding and 
subtracting $M_\pm (\lambda+i0,x_0)$) then yields
\begin{multline}
(M_+(\lambda+i0,x_0)-M_+(\lambda+i0,x_0)^*)
(M_-(\lambda+i0,x_0)-M_+(\lambda+i0,x_0))^{-1} \\
=(M_-(\lambda+i0, x_0)-
M_+(\lambda+i0,x_0))^{-1}(M_-(\lambda+i0,x_0)^*-
M_-(\lambda+i0,x_0)) \lb{4.18}
\end{multline}
and thus
\begin{align}
&(M_-(\lambda+i0,x_0)-M_+(\lambda+i0,x_0))
\Im(M_+(\lambda+i0,x_0)) \lb{4.19} \\
&=-\Im(M_-(\lambda+i0,x_0))
(M_-(\lambda+i0,x_0)-M_+(\lambda+i0,x_0)) 
\text{ for all $\lambda\in\Sigma$}. \no
\end{align}
Taking into account \eqref{4.15} in \eqref{4.19} 
results in
\begin{equation}
(\Im(M_+(\lambda+i0,x_0)))^2=
(\Im(M_-(\lambda+i0,x_0)))^2 \lb{4.20}
\end{equation}
and hence in
\begin{equation}
\Im(M_+(\lambda+i0,x_0))=-\Im(M_-(\lambda+i0,x_0)) 
\text{ for all 
$\lambda\in\Sigma,$} \lb{4.21}
\end{equation}
since $\pm M_\pm(\lambda+i0,x_0)$ are Herglotz 
matrices, implying $\pm\Im(M_\pm (\lambda+i0,x_0))>0.$ 
Combining \eqref{4.15} and
\eqref{4.21}  yields \eqref{4.10} and hence (iii).  

Given (iii) and $M_\pm(z)^*=M_\pm(\bar z)$ one computes 
from \eqref{4.12a}
\begin{equation}
G(\lambda+i0,x,x)^*=-G(\lambda+i0,x,x), \lb{4.22}
\end{equation}
and hence
\begin{equation}
\Xi(\lambda,x)=\f12 I_m 
\text{ for all $x\in\bbR$ and all $\lambda\in\Sigma,$} 
\lb{4.23}
\end{equation}
since $\theta_1(\lambda,x,x_0)$ and 
$\phi_1(\lambda,x,x_0)$ are self-adjoint 
for all $(\lambda,x,x_0)\in\bbR^3.$ Thus (iii) 
implies (i).
\end{proof}

Theorem~\ref{thm4.3} extends to more general situations 
(not necessarily periodic ones) as is clear from the 
corresponding results in \cite{GKT96}, \cite{Ko84}, 
\cite{Ko87a}, \cite{Ko87b}, \cite{KK88}, \cite{SY95} in 
the scalar case $m=1$ 
(replacing the phrase ``for all $\lambda\in\Sigma$'' by 
``for~a.e.~$\lambda\in\Sigma$'', etc.). For the corresponding 
matrix-valued case we refer to \cite{KS88}. 

Combining Theorems~\ref{thm4.8} and \ref{thm4.3}, one 
finally obtains the following result.

\begin{theorem} \lb{thm4.7a}
Suppose Hypothesis~\ref{hyp4.7}. If $H$ has uniform
spectral multiplicity $2m$,
then $H$ is reflectionless and for all $x\in\bbR$ 
and all $\lambda\in\sigma(H)^o,$
\begin{equation}
\Xi(\lambda,x)=\f12 I_m. \lb{4.54}
\end{equation}
\end{theorem}

\begin{corollary} \lb{cor4.9}
Assume Hypothesis~\ref{hyp4.7}. If
$H$ has uniform spectral multiplicity $2m$ and
$\sigma(H)=[E_0,\infty)$ for some $E_0\in\bbR$, then
\begin{equation}
Q(x)=E_0 I_m \text{ for all $x\in\bbR$}. \lb{4.64}
\end{equation}
\end{corollary}
\begin{remark} \lb{rem4.10}
The assumption of uniform (maximal) spectral 
multiplicity $n=2m$ in
Corollary~\ref{cor4.9}
is an essential one. Otherwise, one can easily 
construct nonconstant
potentials $Q(x)$
 such that the associated operator $H$ has overlapping 
band
spectra and hence
spectrum a half-line. For
such a construction it
suffices to consider the case where $Q(x)$ is
a diagonal matrix.
\end{remark}

\begin{remark} \lb{rem4.11}
Corollary~\ref{cor4.9} for matrix-valued 
Schr\"{o}dinger operators
(assuming $Q\in
L^{\infty}(\bbR)^{m \times m}$ to be periodic) has 
been proved by
D\'{e}pres \cite{De95}
using an entirely different approach based on a 
detailed Floquet analysis.
D\'{e}pres' result
partly motivated our work once it became clear that 
trace formulas would be
a most natural tool for
Borg-type uniqueness results and such theorems are 
naturally considered in
the context of
reflectionless rather than periodic potentials. 
Different proofs of Borg's
Theorem~\ref{thm4.1}
(i.e., the scalar case
$m=1$ in Corollary~\ref{cor4.9}) have also been 
obtained in \cite{Jo82},
\cite{Jo87}, \cite{Ko84}.
 Moreover, it
should be stressed that Theorem~\ref{thm4.8} is 
well-known in the special case $m=1$.
(For $m=1$ the hypothesis of uniform spectral multiplicity
$2$ is automatically fulfilled as a simple consequence of
$\det(\Phi(z,x_0))=1$.)  In fact, \eqref{4.55} for 
$m=1$ is proved, for
instance, in  \cite{HS85},
\cite{HKS90}, \cite{Jo82}, \cite{Jo87}, \cite{Ko84}, 
\cite{Ko87a},
\cite{Ko87b}, \cite{SY95}. The reflectionless property
\eqref{4.15} in the general matrix case where
$m\in\bbN$ has also been isolated in \cite{KS88} in 
connection with a study of
stochastic Schr\"{o}dinger and Jacobi operators on strips 
in terms of
Lyapunov exponents.
\end{remark}

\begin{remark} \lb{rem4.12}
In the scalar case $m=1$, the trace formula \eqref{3.21} 
(more precisely, its heat kernel
variant using a heat kernel regularization as opposed 
to our resolvent regularization)
was proved in \cite{GHSZ95} under the general 
condition $V_+ \in
L^1_{\loc}(\bbR)$,
$\sup_{n\in\bbN}(\int_{-n}^n dx\,V_-(x)) < \infty$ 
(where $V_{\pm}=(|V|\pm V)/2$) for all
Lebesgue points $x$ for $V$ (cf.~also \cite{Ry99a} in 
this context). Together with the 
reflectionless property \eqref{4.2} for periodic
potentials this slightly extends Borg's original 
Theorem~\ref{thm4.1} to the case where
$V\in L^1_{\loc}(\bbR)$ is real-valued and periodic. 
We expect a similar
extension of Theorems~\ref{thm3.3} to hold for $x$ in 
the intersection of the Lebesgue points for
$Q_{j,k}$ but did not work out the details.
\end{remark}

At the end we emphasize that all results presented in 
this paper also apply
to matrix-valued
finite-difference Hamiltonian systems.  We refer the 
reader to \cite{CGR98}
in this
direction.

Finally, Borg-type uniqueness theorems for Hamiltonian 
systems are just a
beginning.  There is  a
natural extension of Borg's Theorem~\ref{thm4.1} to 
self-adjoint periodic
Schr\"{o}dinger
operators with one gap in its spectrum, that is,
$\sigma(H)=[E_0,E_1]\cup[E_2,\infty)$, with
$E_1<E_2$.  This extension is due to Hochstadt \cite{Ho65} 
and the
resulting potential $q(x)$
becomes twice the elliptic Weierstrass function. Details 
on a matrix-valued 
extension of Hochstadt's uniqueness theorem will appear 
elsewhere \cite{BGM99}.\\

\noindent {\bf Acknowledgments.} We are very grateful 
to M.~Malamud for a
critical  reading of a preliminary version of this 
manuscript. 
F.~G. and H.~H. gratefully acknowledge financial 
support by the
Research Council of Norway.



\end{document}